\documentclass[a4paper,12pt]{article}
\usepackage[utf8x]{inputenc}
\usepackage{amsmath,amssymb}
\usepackage{a4}
\usepackage{enumitem}

\def\div{\mathop{\rm div}}
\def\RR{\mathbb{R}}
\def\eps{\epsilon}
\def\A{\mathcal{A}}

\newtheorem{theorem}{Theorem}

\begin{document}

\title{Identification of nonlinear heat conduction laws}
\author{H.~Egger\thanks{AG Numerik, Fachbereich Mathematik, TU Darmstadt, Dolivostr. 15, 64293 Darmstadt. Email: {\tt $\{$egger,pietschmann,schlottbom$\}$@mathematik.tu-darmstadt.de}} \and J.~F. Pietschmann\footnotemark[1] \and M.~Schlottbom\footnotemark[1]}
\date{}
\maketitle


\begin{abstract}
We consider the identification of nonlinear heat conduction laws in stationary and instationary heat transfer problems. 
Only a single additional measurement of the temperature on a curve on the boundary is required to determine the unknown parameter function on the range of observed temperatures.
We first present a new proof of Cannon's uniqueness result for the stationary case, then derive a corresponding stability estimate, and finally extend our argument to instationary problems. 
\end{abstract}

{\footnotesize
{\noindent \bf Keywords:} 
parameter identification,
nonlinear diffusion,
quasilinear parabolic problems
}


{\footnotesize
\noindent {\bf AMS Subject Classification:}  
35R30, 65J22
}

\section{Introduction}

This note is concerned with parameter identification problems in nonlinear 
heat transfer processes. 
Let us consider the quasilinear elliptic problem 
\begin{align*}
-\div(a(u) \nabla u) = f \quad \text{in } \Omega, 
\qquad  
a(u) \partial_n u = j  \quad \text{on } \partial\Omega.
\end{align*}
Following \cite{Cannon67}, see also \cite{EggerPietschmannSchlottbom14a,PilantRundell88}, 
the parameter function $a(u)$ can be uniquely determined from temperature measurements $g=u|_{\gamma}$ on a boundary curve $\gamma$.
We present an alternative proof of this uniqueness result below which allows us to treat
also perturbations in the data $f$, $j$, and in $g$, and to obtain a stability result for the 
inverse problem. By using a proper experimental setup, we can then consider also
parabolic problems of the form
\begin{align*}
u_t - \div( a(u) \nabla u) &= f  \qquad \text{on } \Omega \times (0,T),
\end{align*}
with  $u = u_0$ on $\Omega \times \{0\}$ 
and $a(u) \partial_n u = j$ on $\partial\Omega \times (0,T)$. 
In fact, the additional term $u_t$ will be treated as a perturbation of the 
stationary equation.
Our main result about identifiability in the parabolic case can be summarized as follows: 
For any $\eps>0$ and any interval $[g_1,g_2]$, 
we can choose an experimental setup, i.e., data $u_0$, $j$, $f$, and a time horizon $T$, 
such that 
$$
u=\tilde u \text{ on } \gamma \times (0,T) 
\quad  \text{implies that} \quad 
|a(s) - \tilde a(s)| \le \eps \text{ for all } s \in [g_1,g_2].
$$
Here $u,\tilde u$ denote the temperature distributions for parameters $a,\tilde a$, respectively.
It is thus possible to identify the coefficient function $a(u)$ with any prescribed accuracy by a single measurement truely instationary experiment.
Our proof of this result is based on the combination of an energy estimate and a perturbation argument. Similar energy estimates have been used recently also for parameter identification in linear elliptic equations \cite{HarrachSeo2010}.
%
%
%

Parameter identification in heat conduction has a long history \cite{Alifanov94,BeckBlackwellStClair85,Cannon84}. 
To date, rigorous uniqueness results for quasilinear parabolic problems are however only
available in one space dimension \cite{Cannon73,Cannon80,Cannon87,Cortazar90,DuChateau81,DuChateau04};
but see \cite{Isakov93} for multidimensional semilinear problems.
To the best of our knowledge, the question of identifiability of a nonlinear heat conduction law in the multidimensional quasilinear parabolic problem has not been answered yet.
Motivated by applications, several papers are also concerned with numerical methods for parameter estimation in nonlinear heat transfer, see e.g. \cite{Chavent74,Lorenzi86,PilantRundell89}.
For an overview about available uniqueness results and further references on parameter identification in the context of partial differential equations, let us refer to \cite{Isakov06} and \cite{KlibanovTimonov04}.

The remainder of this note is organized as follows: 
In Section~\ref{sec:elliptic}, we introduce some basic assumptions and 
then present our new proof of the uniqueness result of \cite{Cannon67} for the stationary case.
In Section~\ref{sec:stability}, we then derive the corresponding stability result for the inverse problem,
which enables us to treat also quasilinear parabolic equations in Section~\ref{sec:parabolic}. 
We 
conclude with a short discussion. 

\section{Uniqueness for the elliptic problem} \label{sec:elliptic}

Throughout the text, we will make some general assumptions that allow us to keep the presentation simple. Concerning the geometry, we assume that 
\begin{itemize}[leftmargin=2.8em]
 \item[(A1)] $\Omega \subset \RR^d$ is a bounded domain with $C^{1,1}$ boundary in $d=2,3$ space dimensions 
             and $\gamma:[0,1] \to \partial\Omega$ is a $C^1$ curve on the boundary.
\end{itemize}
It should become clear from our analysis that the geometric regularity conditions 
can be further relaxed.
Let us consider the following quasilinear elliptic problem
\begin{align}
-\div (a(u) \nabla u) &= 0 \quad \text{in } \Omega, \label{eq:elliptic1}\\
a(u) \partial_n u &= j \quad \text{on } \partial\Omega. \label{eq:elliptic2} 
\end{align}
In order to ensure the well-posedness of this forward problem, we 
require some regularity and compatibility conditions for the parameter and the data, namely
\begin{itemize}[leftmargin=2.8em]
 \item[(A2)] $a \in \A_{ad} = \{q \in W^{1,\infty}(\RR) : 0 < \underline a
\le q(s) \le \overline{a}$ and $|q'(s)| \le C_1$ for  $s \in \RR\}$. We  denote by $A(s) = \int_0^s a(r) dr$ the principal of $a$.
 \item[(A3)] $f \in L^\infty(\Omega)$, $j \in W^{1,\infty}(\partial\Omega)$, and $\int_\Omega f \; dx + \int_{\partial\Omega} j \; ds(x)= 0$.
\end{itemize}
Via the transformation $U = A(u)$, the quasilinear problem \eqref{eq:elliptic1}--\eqref{eq:elliptic2} can be transformed into a 
Neumann problem for the Poisson equation, and solvability follows from standard results 
for linear elliptic equations \cite{GilbargTrudinger01}; see also \cite{EggerPietschmannSchlottbom14a} for details concerning this particular problem.
We thus obtain
\begin{theorem}
Let (A1)--(A3) hold.
Then \eqref{eq:elliptic1}--\eqref{eq:elliptic2} has a solution $u \in W^{2,p}(\Omega)$ for all $p < \infty$ which is unique up to constants. In addition,  the a-priori estimate $\|u\|_{W^{2,p}(\Omega)} \le C_p (\|j\|_{W^{1,\infty}(\partial\Omega)} + \|u\|_{L^\infty(\gamma)})$
holds with a constant $C_p$ that depends only on $p$, on the geometry, and on the constants of the assumptions.
\end{theorem}
Note that we get uniqueness and a true a-priori estimate once the solution is fixed by the additional temperature measurement $g=u|_\gamma$. It follows from standard embedding results, 
that $u$ and $\nabla u$ are H\"older continuous up to the boundary.

Let us now turn to the parameter estimation problem:
We denote by  $u,\tilde u$ solutions of the elliptic problem \eqref{eq:elliptic1}--\eqref{eq:elliptic2} for parameters $a$, $\tilde a \in \A_{ad}$ with principals $A$, $\tilde A$, and with identical data $f$ and $j$. 
Then for all smooth functions $\phi$, we have 
\begin{align*}
0 
&= -(\div (a(u) \nabla u) - \div(\tilde a(\tilde u)) \nabla \tilde u, \phi )_\Omega \\
&= -(\Delta A(u) - \Delta \tilde A(\tilde u), \phi)_\Omega 
 = (\nabla A(u) - \nabla \tilde A(u), \nabla \phi)_\Omega.
\end{align*}
Here and below, we write $(u,v)_\Omega=\int_\Omega u v \; dx$ for the $L^2$ scalar product. In the last step, we used integration-by-parts and the identical Neumann data 
$\partial_n A(u) = a(u) \partial_n u = j =  \tilde a(\tilde u) \partial_n \tilde u = \partial_n \tilde A(\tilde u)$.
Setting $\phi = A(u) - \tilde A(\tilde u)$ now implies
$\|\nabla A(u) - \nabla \tilde A(\tilde u)\|_{L^2(\Omega)}^2 = 0,$
and hence $A(u) = \tilde A(\tilde u) + c$ with some constant $c \in \RR$. Using the continuity of 
$u$ and $\tilde u$ up to the boundary and assuming identical temperature measurements 
$u|_\gamma = \tilde u|_\gamma=g$, we get 
\begin{align*}
A(g) = \tilde A(g) + c \qquad \text{on } \gamma.
\end{align*}
By differentiation along the curve $\gamma$, we obtain Cannon's uniqueness result \cite{Cannon67}. 
\begin{theorem}
Let (A1) hold and let $u,\tilde u$ denote the solutions of \eqref{eq:elliptic1}--\eqref{eq:elliptic2} for parameters $a,\tilde a \in \A_{ad}$ with identical data $f,j$ satisfying (A3). 
Then the measurement $u|_\gamma=g=\tilde u|_\gamma$ on $\gamma$ implies 
that $a(s) = \tilde a(s)$ for  $s \in 
\text{int} \{g(\gamma(s)) : s \in [0,1]\}$. 
\end{theorem}
Note that the interval of identifiability is empty, if the temperature $g=u|_\gamma$ is constant on $\gamma$. In fact, no identification is possible in that case. 
%

\section{Stability for the inverse problem} \label{sec:stability}

As a second step of our analysis, we now investigate the stability of the identified parameter $a$
with respect to perturbations in the data $f$, $j$, and in the measurements $g$. 
As above, let $u,\tilde u$ denote the solutions of \eqref{eq:elliptic1}--\eqref{eq:elliptic2} for parameters $a,\tilde a \in \A_{ad}$, and with data $f,\tilde f$ and $j,\tilde j$ satisfying assumption (A3). 
Proceeding as in the previous section, we then obtain the identity
\begin{align*}
(f - \tilde f, \phi)_\Omega + (j - \tilde j, \phi)_{\partial\Omega} = (\nabla A(u) - \nabla \tilde A(\tilde u), \nabla \phi)_{\Omega}
\end{align*}
for all smooth functions $\phi$.
Choosing $\phi = A(u) - \tilde A(\tilde u)$ as before and applying the Cauchy-Schwarz inequality to estimate the terms on the left hand side, we get
\begin{align*}
\|\nabla A(u) - \nabla \tilde A(\tilde u)\|_{L^2(\Omega)}^2 
& \le \|f-\tilde f\|_{L^2(\Omega)} \|A(u)-\tilde A(\tilde u)\|_{L^2(\Omega)}  \\
&\qquad \qquad  + \|j - \tilde j\|_{L^2(\partial\Omega)} \|A(u) - \tilde A(\tilde u)\|_{L^2(\partial\Omega)}. 
\end{align*}
Without loss of generality, we can define the principals $A$ and $\tilde A$
in such a way that $\int_\Omega A(u)-\tilde A(\tilde u) dx = 0$. 
By means of the Poincar\'e inequality, we can then deduce that
$\|A(u) - \tilde A(\tilde u)\|_{L^2(\Omega)} \le C_P \|\nabla A(u) - \nabla \tilde A(\tilde u)\|_{L^2(\Omega)}$ 
and similarly that $\|A(u) - \tilde A(\tilde u)\|_{L^2(\partial\Omega)} \le C_P \|\nabla A(u) - \nabla \tilde A(\tilde u)\|_{L^2(\Omega)}$.
Thus we arrive at 
\begin{align*}
\|\nabla A(u) - \nabla \tilde A(\tilde u)\|_{L^2(\Omega)} 
 &\le C_P (\|f-\tilde f\|_{L^2(\Omega)} + \|j-\tilde j\|_{L^2(\partial\Omega)}). 
\end{align*}
Using the uniform boundedness of $A(u)$ and $\tilde A(\tilde u)$ in $W^{1,p}(\Omega)$ with $p$ arbitrarily large, 
interpolation between $L^2(\Omega)$ and $W^{1,p}(\Omega)$ \cite{GriepentrogEtAl02}, embedding of $W^{\theta,q}(\Omega)$ into $C(\overline{\Omega})$ \cite{Adams75}, 
and moving to the boundary, we get
\begin{align*}
\|\nabla A(g) - \nabla \tilde A(\tilde g)\|_{L^\infty(\gamma)} 
\le C_\beta  (\|f-\tilde f\|_{L^2(\Omega)} + \|j-\tilde j\|_{L^2(\partial\Omega)})^{\beta}
\end{align*}
for all $0 \le \beta<3/5$ and $C_\beta$ depending only on $\beta$, on the domain, and the bounds for the coefficients and the data. 
For the choice $\beta=1/2$, we obtain via the triangle inequality, the assumption on the set of addmissible parameters, and by selecting only a tangential component of the gradient
\begin{align*}
&\|\partial_\tau A(g) - \partial_\tau \tilde A(g)\|_{L^\infty(\gamma)} \\
& \qquad \le C_\beta  (\|f-\tilde f\|_{L^2(\Omega)} + \|j-\tilde j\|_{L^2(\partial\Omega)})^{1/2} + \bar a \; \|g-\tilde g\|_{W^{1,\infty}(\gamma)}.
\end{align*}
Here $\partial_\tau$ denotes the derivative along $\gamma$. 
To proceed further, let us assume that 
\begin{itemize}[leftmargin=2.8em]
 \item[(A4)] $[g_1,g_2] \subset \{g(\gamma(s)) : s \in [0,1]\}$ and 
$|\partial_\tau g| \ge c>0$.
\end{itemize}
Note that this is only a technical condition that can easily be satisfied by a proper experimental setup.
By combining the previous estimates, we then obtain
\begin{theorem}
Let $u,\tilde u$ be defined as above and let (A1)--(A4) hold. 
In addition, assume that $\|f-\tilde f\|_{L^2(\Omega)} + \|j-\tilde j\|_{L^2(\partial\Omega)} \le c \eps^{2}$, and $\|g-\tilde g\|_{W^{1,\infty}(\gamma)} \le c \eps$
with some constant $c$ sufficiently small. 
Then  $|a(s) - \tilde a(s)| \le \eps$ for all $s \in [g_1,g_2]$.  
\end{theorem}
The constant $c$ in this theorem only depends on the geometry and the bounds for the coefficients and the data.
If $f=\tilde f$ and $j = \tilde j$, 
we obtain Lipschitz continuity of the parameter $a$ with respect to perturbations in the measurements $g$; 
compare also with the stability result proven in \cite{EggerPietschmannSchlottbom14a}.

\section{Identification in the parabolic case} \label{sec:parabolic}

We will now demonstrate how the argument of the previous section can be extended to parabolic problems of the form 
\begin{align}
 u_t - \div (a(u) \nabla u) &= f \quad \text{on } \Omega \times (0,T), \label{eq:parabolic1}\\
a(u) \partial_n u &= j \quad \text{on } \partial\Omega \times (0,T). \label{eq:parabolic2} 
\end{align}
To ensure the unique solvability and uniform a-priori estimates, 
we assume that
\begin{itemize}[leftmargin=2.8em]
 \item[(A5)] $u(x,0)\!=\!0$, $j \in W^{1,\infty}(\partial\Omega \times (0,T))$ with $j(x,0)\!=\!0$, and $f \in L^{\infty}(\Omega \times (0,T))$.
\end{itemize}
Note that more general initial conditions could be incorporated easily 
and the regularity requirements for the $f$ and $j$  could be further relaxed. 
By standard solvability results for quasilinear parabolic problems \cite{Ladyshenskaja68}, we obtain
\begin{theorem}
Let the assumptions (A1)--(A2) and (A5) hold. 
Then  \eqref{eq:parabolic1}--\eqref{eq:parabolic2} has a unique 
solution $u \in L^2(0,T;W^{2,p}(\Omega))$ for all $p < \infty$
that satisfies
\begin{align*}
&\|u\|_{L^\infty(0,T;W^{1,p}(\Omega))} + \|u\|_{L^2(0,T;W^{2,p}(\Omega))} +\|u_t\|_{L^2(0,T;L^p(\Omega))} \\
&\qquad \qquad \le C_p ( \|f\|_{L^2(0,T;L^\infty(\Omega))} + \|j\|_{W^{1,\infty}(\Omega \times (0,T))}  +  \|u\|_{L^\infty(\gamma \times (0,T))}).
\end{align*}
The constant $C_p$ in the estimate depends only on $p$, on the geometry, and on the bounds for the coefficients and the data used in the assumptions.
\end{theorem}
By  embedding theorems, $u$ and $\nabla u$ are even
H\"older continuous on $\overline{\Omega} \times [0,T]$.

Let us now turn to the parameter estimation problem: 
As before, we denote by $u$ and $\tilde u$ the solutions of \eqref{eq:parabolic1}--\eqref{eq:parabolic2} for parameters 
$a$ and $\tilde a$ with identical data $f$ and $j$. 
Proceeding like in the elliptic case, we obtain for every $0<t \le T$ 
\begin{align*}
(u_t - \tilde u_t, \phi)_\Omega 
&= - (\Delta A(u) - \Delta \tilde A(\tilde u),\phi)_\Omega 
 = (\nabla A(u) - \nabla \tilde A(\tilde u), \nabla \phi)_\Omega.
\end{align*}
Choosing $\phi = A(u) - \tilde A(\tilde u)$ and applying the Cauchy-Schwarz inequality in order to estimate the left hand side, we further get
\begin{align*}
\|u_t - \tilde u_t\|_{L^2(\Omega)} \|A(u)-\tilde A(\tilde u)\|_{L^2(\Omega )}
&\ge  \|\nabla A(u) - \nabla \tilde A(\tilde u) \|_{L^2(\Omega)}^2.
\end{align*}
We can define the principals $A$ and $\tilde A$ in such a way,
maybe differently for every point in time, such that $\int_\Omega A(u)-\tilde A(\tilde u) dx = 0$.
By the Poincar\'e inequality, we then have 
$\|A(u) - \tilde A(u)\|_\Omega \le C_P \|\nabla A(u) - \nabla \tilde A(\tilde u)\|_\Omega$. 
Using this in the previous estimates, we arrive at
\begin{align*}
\|\nabla A(u) - \nabla \tilde A(\tilde u)\|_\Omega \le C_P^{-1} \|u_t - \tilde u_t\|_\Omega 
\qquad \text{for any } 0 < t \le T.
\end{align*}
Note that this inequality is local in time and independent of the choice of the principals $A$ and $\tilde A$. 
By the parabolic nature of the problem, the temperature distribution 
converges exponentially fast to that of the corresponding stationary problem, 
if we keep the data $j$ and $f$ constant over some time and assume that 
they satisfy the compatibility condition (A3).
A slow variation of $f$ and $j$ over time also implies a slow variation of the temperature distribution.
%
By a careful design of the experiment, we may therefore always assume that
\begin{itemize}[leftmargin=2.8em]
 \item[(A6)] $j$, $f$, and $T$ are chosen such that for all $0 < t \le T$ we have $\|u_t\|_{L^2(\Omega)} \le c_1 \eps^{2}$, $|\partial_\tau g| \! \ge \! 1/c_1$,  
  and in addition $[g_1,g_2] \! \subset \! \{u(\gamma(s),t) \! : \! s  \! \in \! [0,1],  t \in [t_1,t_2]\}$ for some $0 < t_1 \le t_2 \le T$. 
\end{itemize}
Using this condition, the uniform a-priori estimates for the solution, 
an interpolation argument, and moving to the boundary, we conclude
\begin{theorem}
Let (A1) and (A5) hold and denote by $u,\tilde u$ the solutions of problem \eqref{eq:parabolic1}--\eqref{eq:parabolic2} with parameters $a,\tilde a$ satisfying (A2). 
Moerover, assume that the experiment is designed such that (A6) is valid with $c_1$, $c_2$ sufficiently small.
Then from measurements $u(x,t) = \tilde u(x,t) = g(x,t)$ on $\gamma \times [t_1,t_2]$, 
we may conclude that $|a(s) - \tilde a(s)| \le \eps$ for all $s \in [g_1,g_2]$.
\end{theorem}
Note that the constants $c_1$, $c_2$ do only depend on the domain and the bounds for the coefficients and the data. 
It should become clear from the derivation above 
that the statement of the theorem can be localized in time, i.e., we may 
identify $a(s)$ on $[g_1,g_2]$ from measurements $g=u(\gamma,t^*)$ at a single point in time 
on the corresponding range of temperatures. 
Similar as in the elliptic case, 
perturbations in the measurements and the data could be incorporated as well.

\section{Discussion}

In this note, we investigated the identification of nonlinear heat conduction laws in 
stationary and instationary heat transfer processes. We presented a new proof for Cannon's uniqueness result for the quasilinear elliptic problem which allowed us to derive 
a corresponding stability result with respect to perturbations in the data. 
Using this stability result, the parabolic problem could then be treated as a perturbation of the elliptic case. We finally could obtain the approximate
identification of the unknown parameter function with arbitrary accuracy 
by a single experiment.
Let us mention that, in principle, one could also apply a time-independent temperature flux, wait until the system reaches equilibrium, and then apply the results for the perturbed elliptic problem. By the parabolic nature, the system will reach the stationary equilibrium exponentially fast. In contrast to such an approach, the setting considered in Section~\ref{sec:parabolic} is truely instationary, but close to the stationary equilibrium for all times. 

\section*{Acknowledgments}
The first author acknowledges support by DFG via grants IRTG~1529 and GSC~233.
The work of the second author was supported by DFG via grant 1073/1-1 and by the Daimler~and~Benz Stiftung via stipend 32-09/12. 
%

\end{document}